\def\printVersionNumber{\rightline{version 2.5}}

%
\def\printVersionNumber{}




 \message{Assuming 8.5 x 11 inch paper.}
 \message{' ' ' ' ' ' ' '}

\magnification=\magstep1	          
\raggedbottom

%

\parskip=9pt

\def\singlespace{\baselineskip=12pt}      
\def\sesquispace{\baselineskip=15pt}      


\font\titlefont=cmb12 at 12pt



%



\font\openface=msbm10 at10pt
 %

\def\Complexes     {{\hbox{\openface C}}}
\def\Minkowski     {{\hbox{\openface M}}}

 %
 %
 %



\font\german=eufm10 at 10pt

\def\Buchstabe#1{{\hbox{\german #1}}}





\def\Aut {\mathop{\rm Aut} \nolimits}	 


\def\span{\mathop {\rm span }\nolimits}

%

%
%



\def\implies{\Rightarrow}

%


 \def\dal{\displaystyle{{\hbox to 0pt{$\sqcup$\hss}}\sqcap}}



\def\lto{\mathop
        {\hbox{${\lower3.8pt\hbox{$<$}}\atop{\raise0.2pt\hbox{$\sim$}}$}}}
\def\gto{\mathop
        {\hbox{${\lower3.8pt\hbox{$>$}}\atop{\raise0.2pt\hbox{$\sim$}}$}}}
%
%
%



\def\&{{\phantom a}}

\def\braces#1{ \{ #1 \} }

\def\bra{<\!}			
\def\ket{\!>}			

\def\isom{\simeq}		

\def\to{\mathop\rightarrow}	


\def\SetOf#1#2{\left\{ #1  \,|\, #2 \right\} }


\def\union{\cup}

\def\cross{\times}

\def\interior #1 {  \buildrel\circ\over  #1}     




\def\basisvector#1#2#3{
 \lower6pt\hbox{
  ${\buildrel{\displaystyle #1}\over{\scriptscriptstyle(#2)}}$}^#3}





%


%
%
%
%
%
%
%
%
%
%

%
 \let\miguu=\footnote
 \def\footnote#1#2{{$\,$\parindent=9pt\baselineskip=13pt%
 \miguu{#1}{#2\vskip -7truept}}}
%
%

\def\linebreak{\hfil\break}


\def\BulletItem #1 {\item{$\bullet$}{#1 }}

\def\AbstractBegins
{
 \singlespace                                        
 \bigskip\leftskip=1.5truecm\rightskip=1.5truecm     
 \centerline{\bf Abstract}
 \smallskip
 \noindent	
 } 
\def\AbstractEnds
{
 \bigskip\leftskip=0truecm\rightskip=0truecm       
 }

\def\ReferencesBegin
{
 \singlespace					   
 \vskip 0.5truein
 \centerline           {\bf References}
 \par\nobreak
 \medskip
 \noindent
 \parindent=2pt
 \parskip=6pt			
 }

\def\section #1 {\bigskip\noindent{\bf #1 }\par\nobreak\smallskip\noindent}

\def\subsection #1 {\medskip\noindent[ {\it #1} ]\par\nobreak\smallskip}

\def\reference{\hangindent=1pc\hangafter=1} 

\def\ref{\reference}

\def\journaldata#1#2#3#4{{\it #1 } {\bf #2:} #3 (#4)}
 %

\def\eprint#1{$\langle$#1\hbox{$\rangle$}}
 %


\def\author#1 {\medskip\centerline{\it #1}\smallskip}

\def\address#1{\centerline{\it #1}\smallskip}

\font\titlefont=cmb10 scaled\magstep2 



\def\THEOREM{\noindent {\csmc Theorem \ }}
\def\LEMMA{\noindent {\csmc Lemma \ }}





\def\A{{\Buchstabe A}}			
\def\B{{\Buchstabe B}}			
\def\I{{\Buchstabe I}}			

\def\U{{\cal U}}			

\def\LEMMA{\noindent LEMMA \ }
\def\PROOF{\noindent PROOF \ }
\def\THEOREM{\noindent THEOREM \ }
\def\COROLLARY{\noindent COROLLARY \ }
\def\DEFINITION{\noindent DEFINITION \ }

\def\up{upward}		
\def\down{downward}

\def\Inner{\mathop {{\rm Inner}} \nolimits}
\def\Outer{\mathop {{\rm Outer}} \nolimits}


\phantom{}
\vskip -1 true in
\medskip

\printVersionNumber

\rightline{gr-qc/0309126}
\rightline{SU--GP--03/6--3}	   

\vskip 0.3 true in
\bigskip
\bigskip

\sesquispace
\centerline{ {\titlefont 
  Indecomposable Ideals in Incidence Algebras}\footnote{$^\star$}%
{To appear in a special issue of {\it Modern Physics Letters A}
 devoted to the proceedings of ``Balfest'',
 held May, 2003, in Vietri sul Mare, Italy}}

\bigskip


\singlespace			        

\author {Rafael D. Sorkin}
\address
 {Department of Physics, Syracuse University, Syracuse, NY 13244-1130, U.S.A.}
\centerline {\it \qquad\qquad internet email address: sorkin@physics.syr.edu}

\AbstractBegins
The elements of a finite partial order $P$ can be identified with the
maximal indecomposable two-sided ideals of its incidence algebra $\A$,
and then for two such ideals, $I\prec J \iff IJ \not=0$.  This offers
one way to recover a poset from its incidence algebra.
In the course of proving the above, we classify all of the two-sided
ideals of $\A$.
\AbstractEnds


\sesquispace


\noindent
In contemporary physical theory, the concept of a ``space'' or, more
formally, of a set with structure, plays a central role.  Most notably,
spacetime itself is conceived of in this manner -- as a differentiable
manifold.  However, one can observe a certain tension between two ways
of conceptualizing such structures and working with them.  Tangent
vectors, for example, can be 
thought of either as infinitesimal
displacements or as sets of numbers obeying a linear
transformation law.  In present day language, the two opposed 
tendencies 
of thought
can
to some extent be characterized by the words ``geometrical'' and
``algebraic'', although neither term is really suitable.  Perhaps,
``intrinsic'' vs ``coordinate based'' comes closer; 
and sometimes the words ``synthetic'' and ``analytic'' have been used to
convey the same opposition (as in synthetic versus analytic
geometry).\footnote{$^\dagger$}
{In the philosophical world, these two attitudes manifest themselves to
 some degree as ``materialism'' and ``instrumentalism'' although the
 correspondence is obviously very imperfect (cf. the oft stated idea that
 an instrument reading is always a number.)}

Consider, for example, Minkowski spacetime $\Minkowski^4$.  From the
``intrinsic'' side it can be understood, on one hand, as a topological
space of dimension four supporting such concepts as straight line
(inertial motion), light cone, and parallelogram.  Or by focusing on its
causal relationships rather than its metric and topological ones, one
can understand $\Minkowski^4$ as a certain partially ordered set
(poset), a second ``intrinsic'' characterization that is nevertheless
very different from the first.  In contrast to both these
characterizations, $\Minkowski^4$ would be described from the
``coordinate based'' side by four real variables $t,x,y,z$ which
geometrically have the meaning of numerical functions on spacetime.  Here,
the algebraic relationships among the four variables take center stage,
while the actual elements of the space (the points) withdraw into the
background.

Of course the ``intrinsic'' and ``coordinate based'' descriptions of
$\Minkowski^4$ are mathematically equivalent.  The most highly developed
and general instance of this sort of equivalence is the
Gel'fand isomorphism, which implies in particular
that any manifold can be recovered,
as a topological space, 
from the $C^*$-algebra of scalar functions that
it supports (in effect its coordinate functions).  However, a manifold
{\it per se} is not yet a spacetime because it lacks metrical information.  To
recover that as well, one can proceed as in [1] or following
the more detailed scheme of [2].  Neither approach captures in
any essential manner the {\it Lorentzian} character of the metric,
however.  (Indeed, the latter scheme is actually incompatible with
Lorentzian signature.)  The question thus arises whether there exists a
similar ``algebraization'' of spacetime based not on its topological and
metrical attributes but on its causal order.

The finding of such a correspondence could be expected to hold interest
for more than one reason.  On one hand, some workers, going back to
[3], have viewed algebraization as potentially a means by which
to introduce a fundamental spacetime discreteness, a view which appears
to account for much of the current interest in ``non-commutative
geometry''.  On the other hand, algebraization has from the outset been
one of the royal roads to the ``quantization'' of a theory, so that one
might hope that any new equivalence between intrinsic and algebraic
descriptions of spacetime -- or of whatever hypothetical substratum one
takes to replace spacetime -- would open up new avenues for building a
theory of quantum gravity.  It is this second prospect that primarily
animates the considerations of the present paper, which are inspired by
the hypothesis that the deep structure of spacetime is that of a {\it
causal set} [4] [5].  Since this structure is already inherently
discrete, there is no need to {\it introduce} discreteness and therefore
no reason to appeal to algebraization on that score.  Nonetheless, one
may still feel it useful to attempt various algebraic reformulations of
causal set kinematics in the hope that one of them might help lead us to
the correct quantum theory of causal set dynamics.
%
%
I will have a bit more to say about this in the conclusion, but for now
let us turn to the mathematical question that will primarily concern us
herein: that of finding a suitable ``algebraization'' of the poset
concept.  

As I remarked earlier, a relativistic spacetime is inherently a
partial order\footnote{$^\flat$}
{{\it Order}, 
 {\it partial order}, 
 {\it partially ordered set}, 
 {\it poset}, 
 and 
 {\it ordered set} 
 are all synonyms.}
(at least to the extent that one can count on the impossibility of
 ``time travel'').
A causal set is also a partial order, but with the crucial difference
that it is {\it locally finite}.\footnote{$^\star$}
{A poset is called locally finite if all its {\it intervals} are
 finite.  If this is strengthened to the requirement that its ``{\up}''
 and ``{\down}'' sets be separately finite, the poset becomes suitable
 to represent a ``finitary topological space''.  Posets of this type
 have been the subject of much work by the dedicatee and his
 co-workers, but in a spatial context rather than a spacetime one
 [6]. (See also [7] [8]).  
 The results to be proven below will of course apply
 equally well in this context.}
Rather than 
seeking an algebra to capture the structure of an arbitrary poset, then,
let us confine ourselves to the simpler case of
locally finite orders; in fact let us simplify still further
to the case of orders whose cardinality is strictly finite.

So let $C$ be a poset of finite cardinality.
Does there exist an algebra $\A$ naturally associated to $C$ and from
which, conversely, $C$ can be recovered?  One algebra that people have
studied in this connection is the so called {\it incidence algebra} of
$C$, which one might view as an algebra of retarded (or, dually,
advanced) functions on $C{\cross}C$, the product being given by
convolution.  In the finite case, this is just a matrix algebra, or more
accurately, it is the algebra of all matrices with zeros in certain
specified locations which reflect directly the defining order relation
$\prec$ of $C$.

It is known that the incidence algebra does indeed capture the structure of
$C$, at least if one interprets $\prec$ reflexively in the
sense that the diagonals of the matrices representing the elements of
$\A$ are left free.  Remarkably, the nature of the 
space $\leftrightarrow$ algebra correspondence 
in this case 
can be arranged to be 
closely analogous to that of the Gel'fand isomorphism,
even though the algebra which figures 
in the latter 
is commutative and semisimple,
whereas an 
incidence algebra is not only non-commutative, 
but almost nilpotent 
(which is about as far from semi-simplicity as one can get).  
Despite these differences, one can in both cases 
choose to
identify the
elements of the underlying space with the maximal two-sided ideals of
$\A$, as explained in [9] and [10].

As we will see in a moment, however, there is quite a different way to
set up a correspondence between $C$ and $\A$ in the poset case, and this
is perhaps fortunate, in that one might feel uncomfortable with certain
features of the Gel'fand-like scheme as adapted to partial orders.  A
first concern arises from the circumstance that it is not the 
fundamental order
relation $\prec$ as such that one directly recovers between elements
regarded as maximal ideals, but rather the ``nearest neighbor'' or {\it
link} relation\footnote{$^\dagger$}
{Also called ``covering relation'' in the mathematical literature.}
[9] [11].
The full precedence relation must then be re-generated via transitive 
closure.  
In the purely discrete case, this is always possible, so viewing one
relation as more or less ``fundamental'' than the other seems largely a
matter of taste.  
However, the failure to recover $\prec$ directly could bode ill in a
continuum context where the simple precedence relation continues to make
sense but links no longer exist (since between any two causally related
points, one can always interpolate a third).
%
%
Moreover, even for finite posets $C$, the correspondence between
elements and maximal ideals falls apart if one adopts the irreflexive
convention\footnote{$^\flat$}{defined below} for $C$, and this is
disturbing because one would hope that a choice of convention would not
influence the underlying relationships so strongly.  Again, one might
hope that some other scheme would be more robust in this regard.

In view of such doubts, it seems worth exploring other ways to recover a
poset from its incidence algebra, a task we begin here by showing that
it is equally possible to identify the poset elements, not with the {\it
maximal} ideals, but with certain {\it indecomposable} ideals; and 
by doing so
one obtains the precedence relation directly as a relation
between the corresponding ideals.  Of course, it would also make sense
to explore alternatives to the incidence algebra itself, but I will not
attempt that in this paper.

\section{ Some Definitions }

A {\it finite order} is a set $C$ comprising a finite number of
elements and carrying a ``precedence relation'' $\prec$ which is
transitive and asymmetric.  That is, 
for arbitrary elements $x,y,z$ of $C$ 
we always have $x\prec y\prec z\implies x\prec z$, 
and we never have $x \prec y \prec x$ when $x$ and $y$ are  distinct.  
In addition we will
always assume, unless stated otherwise, that $\prec$ is {\it reflexive},
i.e. that every element $x$ precedes itself: $x\prec x$.  Although, this
is in some sense merely a convention that one can adopt or reject at
will, it turns out to influence profoundly the structure of the
incidence algebra.\footnote{$^\star$}
{For example, with the reflexive convention, $\A$ has in general more
 idempotents than $C$ has elements, whereas with the opposite,
 {\it irreflexive} convention, $\A$ has no idempotents at all.}

We can now define the {\it incidence algebra} of $C$ by introducing for
each related pair $x{\prec}y$ a generator $[xy]$ and taking $\A$ to be
the set of formal linear combinations of these generators.  The order 
structure of $C$ is then further encoded in the rule for
multiplication of algebra elements as given by the relation
$[xy] \cdot [yz]=[xz]$, for all triples $x \prec y \prec z$.
For definiteness, we will take the field of scalars to be the complexes
$\Complexes$, although nothing will depend on this.

Two other representations of the incidence algebra are useful 
and (for finite $C$) strictly equivalent to the definition just given.  
First one may think of $\A$ as an algebra of $n\times n$ {\it matrices}, 
where $n=|C|$ is the cardinality of $C$ and the generator $[xy]$ corresponds
to the matrix with a single $1$ in row $x$, column $y$ and 
zeros everywhere else.\footnote{$^\dagger$}
{In the Dirac notation, this correspondence reads 
     $[xy] = |{x}\ket \bra{y}|$, 
 a notation that was used in [9].} 
The asymmetry of the relation $\prec$
then translates into the fact that
if one chooses a suitable labeling for the elements of $C$
(a so called {\it natural labeling})
then
all the matrices representing members of $\A$ are upper triangular 
(but not strictly so, inasmuch as diagonal generators like $[xx]$ are
 also part of $\A$ thanks to our standing assumption that $\prec$ is
 reflexive).
A slightly different
representation of $\A$ comes from
thinking of the entries of the matrix as the values of a ``two-point
function'' $f:C\times C\to\Complexes$.  With this representation a member
of the incidence algebra is an arbitrary {\it advanced} function and the
algebra product is convolution: $f=g\cdot h\iff f(x,z)=\sum_y g(x,y) h(y,z)$. 
Of course this is just the formula for matrix multiplication in a slightly
different notation.

By an {\it ideal} $I$ of $\A$, I will always mean, in this paper, a
two-sided ideal, in other words a nonempty
subset of $\A$ closed under addition,
scalar multiplication,
and left or right multiplication by an arbitrary element of $\A$.  
Note in this connection that, by virtue of our choice of reflexive
convention for $C$, $\A$ is automatically ``unital'': it has an identity
element given by $1=[xx]+[yy]+[zz]+\cdots$, where the sum extends over
all elements of $C$.  Consequently there is no need to distinguish, for
example, ``regular ideals'' from irregular ones [12].

The {\it sum} $I_1+I_2$ of two ideals $I_1$ and $I_2$ is the
collection of all sums $a_1+a_2$ where $a_1\in I_1$ and $a_2\in I_2$.
Equivalently $I_1+I_2$ is the least ideal containing both $I_1$ and
$I_2$: it is their ``join'' in the lattice of ideals of $\A$.  
An {\it indecomposable ideal} 
in $\A$ is then a 
non-zero 
ideal that cannot be
expressed as the sum of two ideals distinct from itself. 
(In the language of lattice theory, such an ideal is said to be ``join
irreducible''.  A closely related notion was employed in the definitions
of $TIP$ and $TIF$ in [13].)  
By a {\it maximal indecomposable ideal} 
I will mean an indecomposable ideal
that is contained properly in no other indecomposable ideal.  
One also defines simply a {\it maximal ideal} $I\not=\A$ as one which
cannot be enlarged without coinciding with $\A$.  
(Thus, a maximal ideal is
a maximal element in the family of all ideals not equal to $\A$,
while a maximal indecomposable ideal is
a maximal element in the family of all indecomposable ideals.)  

Similarly, we define the {\it product} of two ideals $I$ and $J$ 
as the collection of products of members of $I$ with members of $J$:
$IJ=\SetOf{xy}{x,y\in\A}$.  From the definitions it is immediate that
$IJ$ is also an ideal and is contained in both $I$ and $J$.  (It is,
however, not in general equal to their intersection $I\cap J$, which
would be their ``meet'' in the lattice of ideals.)

Finally, we define within an arbitrary order $P$ a {\it {\down} set} as a
subset $D\subseteq P$ that is closed under ``taking of pasts'':
$x\prec y\in D\implies x\in D$.  
An {\it {\up} set} $U$ is defined dually.  
(In spacetime language, these could be called, respectively, ``past
sets'' and ``future sets''.)

\section{ Recovery of the poset from its incidence algebra }

Our main result can be stated as a theorem:

\THEOREM Let $C$ be a finite order (with reflexive convention) and let
$\A$ be its incidence algebra.  Then the elements $x$ of $C$ correspond
bijectively with the maximal indecomposable ideals $I$ of $\A$, and
under this correspondence the relation  $x_1\prec{x_2}$ goes over to
$I_1 I_2 \not= 0$.

Before proving the theorem, let us notice a sense in which this (or any
other) equivalence between a poset and its incidence algebra is somewhat
less satisfactory than the Gel'fand isomorphism mentioned earlier, the
flaw being that $\A$ has in general more automorphisms than $C$
has.\footnote{$^\flat$}
{Because the algebra $\A$ is not commutative, it will in general possess
 continuous families of inner automorphisms, whereas $C$, being a finite
 set, can have at most a finite number of automorphisms.  (On the other
 hand, it looks as if the superfluity due to the inner automorphisms
 might be the only one.  That is, it looks as if we might have
 $\Aut(C)\isom\Outer(\A):=(\Aut\A)/(\Inner\A)$, where $(\Inner\A)$ is
 the group of inner automorphisms of $\A$.)}
Although the isomorphism equivalence class of $\A$ can be deduced from
that of $C$ and conversely, there is thus some ``looseness'' in the
correspondence that is not present in the Gel'fand case.  
%
%
Possibly related is the failure, pointed out in [14], of
the correspondence between a poset and its incidence algebra to be
functorial between the category of finite orders with isotone mappings
and the category of algebras with algebra homomorphisms.\footnote{$^\star$}
{This is not necessarily a ``failing'' in itself.  For example, the
 association to a manifold $M$ of its {\it tensor} algebra is not
 functorial, nor (as remarked to me by Chris Isham) is the association
 to $M$ of its diffeomorphism group.  Nevertheless, for certain purposes
 the lack of functoriality can be a problem, e.g. if one wished to
 reproduce the limiting process of [7] in terms of incidence
 algebras.}
%
%

\section {The ideals of $\A$}
We will prove the above theorem by classifying the ideals of $\A$.  To
this end, let us notice that every member of $\A$ has a unique
expression as a sum of multiples of the generators $[xy]$
that we defined
earlier.  To convey the fact that one particular such generator $[xy]$
occurs with a nonzero coefficient in some member of the ideal $I$,
let me say, for lack of a better word, that $[xy]$ ``figures in $I$''.  
In
contrast, the statement that $[xy]$ ``is an element of $I$''
(in symbols, $[xy]\in I$)
means that
some $A\in{I}$ literally coincides with $[xy]$.  In the matrix
representation of $\A$, a generator $[xy]$ corresponds to a particular
location in the matrix.  It then ``figures in'' $I$ if some matrix of
$I$ has a nonzero coefficient in that location, whereas it is an element
of $I$  if some matrix of $I$ has a $1$ in that
location and zeros everywhere else.  In these terms we can now state a
key lemma.

\LEMMA If some $[xy]$ ``figures in'' the ideal $I$ then it is actually
an element of $I$

\PROOF By assumption there is $A\in I$ such that $A=\alpha[xy]+B$ where
$\alpha$ is a nonzero scalar and 
$B$ is a sum of multiples of pairs $[uv]$ such that either $u$ differs from
$x$ or $v$ from $y$.  But this means that $[xx]B[yy]=0$ because the same
holds for every one of its constituent 
pairs
$[uv]$.  
Hence
$[xx]A[yy]=\alpha[xx][xy][yy]=\alpha[xy]$ 
is a multiple of $[xy]$, and this implies
immediately that $[xy]$ itself belongs to $I$ by the definition of an
ideal.  
(Namely $[xy]= (1/\alpha) [xx]A[yy]$ must be an element of $I$ if
$A$ is one.)

\COROLLARY 
Every ideal $I$ of $\A$ is the set of all linear combinations of 
some unique set 
$\U(I)$
of generators $[xy]$.

\noindent
That is, we have for every ideal of the incidence algebra,
$I=\span\U(I)$, 
where
$\U(I)=\braces{[x_1y_1], [x_2y_2], [x_3y_3], \cdots [x_ny_n]}$, 
with the
$[x_jy_j]$ 
being uniquely determined by $I$ and conversely.  
It is moreover, easy to figure out which 
sets  
of pairs can belong to an
ideal in this way.  
Let $S=\U(I)$ be such a set.  If
$[xy]\in S$ and $u\prec x\prec y\prec v$ then $[uv]=[ux][xy][yv]\in{I}$
and therefore $[uv]\in{S}$.  Thus $S$ is necessarily closed under the
process of ``passing to nested pairs''.  To express this succinctly let
us formally introduce this ``nesting'' as an order relation among pairs.

\DEFINITION $[xy]\ll[uv] \iff u\prec x$ and $y\prec v$

\noindent
This definition makes the set of all pairs
$[xy]$ into a poset $\Gamma$,
and with reference to this auxiliary poset, we see that the sets 
$\U(I)$
are precisely the {\it {\up} sets}
of $\Gamma$.  Thus we have proved:

\THEOREM The ideals of $\A$ correspond bijectively with the {\up} sets of
the poset $\Gamma$ of pairs $[x,y]$

\noindent
Moreover, 
the relation of inclusion between ideals obviously mirrors exactly 
the relation of inclusion between the corresponding {\up} sets.  
In particular, we have

\LEMMA $\U(I+J) = \U(I) \union \U(J)$

\PROOF Recall that $I+J$ is precisely the smallest ideal including both
$I$ and $J$; and notice that the union of two {\up} sets is also an {\up}
set. 

It follows immediately that
an ideal $I$ is indecomposable iff its corresponding {\up} set, $\U(I)$,
is not the union of two {\up} sets distinct from itself.
But in any (finite) poset 
an {\up} set has this property iff it is 
(empty or)
what might be called a
``principal {\up} set'', that is iff it has the form
$\SetOf{\eta}{\xi\ll\eta}$ for some $\xi$.  We have thus shown:

\THEOREM
The {\it indecomposable ideals} of $\A$ correspond bijectively 
with the pairs $[xy]$, $x,y\in C$, $x\prec{y}$


From this we can easily conclude that the 
{\it maximal} indecomposable ideals correspond precisely with the
principal {\up} sets of {\it minimal} elements of $\Gamma$, 
which in turn are clearly 
the diagonal pairs $[xx]$, for $x{\in}C$.  
In other words:

\LEMMA An ideal $I\subseteq\A$ is {\it maximal indecomposable} 
$\iff I = \A [xx] \A$ for some $x\in C$.

This substantiates the first assertion of our main theorem.\footnote{$^\dagger$}
{Very similar reasoning yields the ``Gel'fand'' correspondence between
 maximal ideals and elements of $C$, since a maximal {\up} set in $\Gamma$ is
 precisely the complement of a single minimal element of $\Gamma$, i.e.  (as
 we have just observed) of a single pair $[xx]$.  However, if our only
 purpose were to classify the maximal ideals of $\A$, it would be
 simpler just to adopt the representation of an element of $\A$ as a
 matrix and pay attention to the set of zeros on its diagonal.}
%
%
In order to complete the proof, we have only to verify
that $I_x I_y \not= 0$ iff $x{\prec}y$,
where I've written $I_x$ for $\A[xx]\A$.
But this is completely
straightforward. In fact, we have:

\LEMMA Let $\I(x,y)=\A[xy]\A$ denote the principal ideal generated by $[xy]$.
Then the product $\I(x,y) \I(u,v)$ equals $\I(x,v)$ when $y\prec u$, and
zero otherwise. 

\PROOF Write $\B$ for the product in question.  Because $\A$ is unital,
$\A\A=\A$, whence $\B=\A[xy]\A\A[uv]\A=\A[xy]\A[uv]\A$.  Now, 
$[xy]\A[uv]\not=0 \iff [yu]\in\A \iff y\prec u$, and in 
that
case,
clearly, $\B=  \A[xy][yu][uv]\A   =\A[xv]\A$.

\COROLLARY 
$\I(x,x)\I(y,y)$ equals $\I(x,y)$ when $x{\prec}y$, and zero otherwise.

Our main theorem is thus demonstrated.  Along the way, we have seen that
every indecomposable ideal is a principal ideal (of some pair $[xy]$).
We have also found (in view of the most recent lemma) an equivalent way to
characterize the maximal indecomposable ideals: An indecomposable ideal
$I$ is maximal iff it is ``idempotent'' in the sense that $II=I$
(which happens iff $II\not=0$).
%

Finally, it bears remarking that the last lemma actually furnishes an
order on the full set of indecomposable ideals.  Since every such ideal
has the form $\A[xy]\A$
for some $[xy]$, 
the substrate of this order can be taken to be the set of pairs\footnote{$^\flat$}
{Thus, the order in question shares its substrate with the order
 $\Gamma$.  The two precedence relations are obviously very different,
 however.}
$x\prec y$ of elements of $C$ and its precedence relation is
then
$[xy]\prec[uv]\iff y\prec u$. 
(Such an order is often called an ``interval order''.)  
It is actually this
order that we recover most directly from the algebra $\A$.
Interestingly enough, it is neither reflexive nor
irreflexive.\footnote{$^\star$}
{providing a counter-argument to the opinion that the choice between the
 two possibilities is always purely conventional!}
Rather, it is precisely its reflexive elements that correspond to
elements of the underlying poset $C$.

\section{ Remarks }

With the proof of our main theorem, we are in the curious position of
being able to discern two very different algebraic images of our
original poset $C$ within its incidence algebra $\A$.  On one hand, we
can identify the element $x\in C$ with the ideal {\it generated by}
$[xx]$, on the other hand with the ideal of all algebra elements {\it
omitting} $[xx]$.  Depending on which possibility one selects, one
recovers from the ideals either the relation $\prec$ or its associated
link relation, respectively.  Either way, one can conclude that the full
structure of $C$ is captured by $\A$.

However, with the second scheme, this conclusion rests heavily on our
assumption that $C$ is finite, which (trivially) guarantees enough
discreteness so that any two related elements of $C$ can be joined by a
chain of links.  In a continuum such as a Lorentzian manifold (or for a
countable dense set thereof), this is assuredly not the case because no
links are present at all.  Hence the construction which recovers $\prec$
directly (that based on maximal indecomposable ideals) seems more robust
from this point of view.  It might be interesting to test whether this
is indeed true by generalizing $\A$ to, say, a globally hyperbolic
spacetime $M$ and asking whether some suitable analog of the
construction of this paper would give back the metric of $M$.
(This would require one to recover the volume density $\sqrt{-g}$
in addition to the causal order of $M$, but that is not obviously
impossible, since $\sqrt{-g}$ does enter into $\A$, via the
definition of the convolution product.)

Of course, it would also be interesting to investigate possible
inter-relationships between the two schemes for recovering the elements
of $C$.  Perhaps they could already be seen to be equivalent at the
algebraic level, or perhaps they complement each other in some way.

Concerning physical applications of the duality between a poset and its
incidence algebra, there seems not much to say at present.  If one is
thinking in terms of quantum gravity and causal sets, then trading a
causal set $C$ for its incidence algebra $\A$ is not necessarily a step
in the right direction, because a quantal ``sum over causal sets'' seems
easier to imagine (as in [5]) than a corresponding ``sum over
algebras $\A$''.  On the other hand, the physical causal set $C$ is not
the only poset to figure in the theory.  The ``poset of stems'' played
an important role in the considerations of [5], where it served as
a fixed arena allowing one to define conveniently the Markov process of
``classical sequential growth''.  (This poset is illustrated in Figure 1
of [5].  Its elements are the finite orders, and its precedence
relation is ``inclusion as a {\down} set''.)  Perhaps the incidence
algebra of {\it this} poset could play a role in the search for a
physically appropriate dynamics of ``quantal sequential growth''.  Such
a dynamics would provide the quantal analog for causal sets of the
classical Einstein equations for continuum Lorentzian manifolds.  That
is, it would provide a theory of quantum gravity.

Returning to the realm of pure mathematics for a moment, one can ask
whether algebras other than the incidence algebra might have a role to
play in the ``algebraization'' of the poset concept.  If so, the new
algebra might be constructed either from the incidence algebra or
directly from the elements of the poset.  Here, just to illustrate the
type of thing one could consider, is an algebra of the second sort.  It
might even be trivial as far as I know, but at least it is not as
obviously trivial as some similar possibilities I played with first!  One
takes for generators the elements of $C$ itself, and one imposes three
sets of relations: 
(i) $a\prec b \prec c \implies abc=ac$;
(ii) $a \not= b \implies aba=0$;
(iii) $a \not= b, \, a \prec b \implies ba=0$.
If this construction is worthy of further consideration, one might begin
by asking whether it is functorial and how the resulting algebra's
automorphism group relates to that of $C$.


\bigskip
\noindent
I would like to thank Prakash Panangaden, Ioannis Raptis and Kamesh Wali
for inspiration and encouragement. 
This research was partly supported by NSF grant PHY-0098488 and by a
grant from the Office of Research and Computing of Syracuse University.

\ReferencesBegin

\ref 
[1]
R.~Geroch, 
 ``Einstein Algebras'', 
  \journaldata {Comm. Math. Phys.} {26} {271-275} {1972}

\ref 
[2] 
Alain Connes, {\it Noncommutative Geometry} (Academic Press, 1994);
\linebreak
Giovanni Landi, {\it An Introduction to Noncommutative Spaces and their
 Geometries}, 
 Lecture Notes in Physics: Monographs, m51
 (Springer-Verlag, 1977)
 \eprint{hep-th/9701078}.

\ref 
[3]
Hartland S. Snyder, ``Quantized Space-Time'', 
\journaldata{Phys. Rev.}{71}{38-41}{1947}

\ref 
[4]
L.~Bombelli, J.~Lee, D.~Meyer and R.D.~Sorkin, 
``Spacetime as a causal set'', 
  \journaldata {Phys. Rev. Lett.} {59} {521-524} {1987}

\ref 
[5]
David P.~Rideout and Rafael D.~Sorkin,
``A Classical Sequential Growth Dynamics for Causal Sets'',
 \journaldata{Phys. Rev. D}{61}{024002}{2000}
 \eprint{gr-qc/9904062}

\ref		
[6]
A.P.~Balachandran, G.~Bimonte, E.~Ercolessi, G.~Landi, F.~Lizzi,
   G.~Sparano, and P.~Teotonio-Sobrinho,
``Noncommutative Lattices as finite approximations'',
   \journaldata {\it J. Geom. Phys.} {18} {163} {1996}
   \eprint{hep-th/9510217}

\ref 
[7]
R.D.~Sorkin, ``A Finitary Substitute for Continuous Topology?'',
  \journaldata {Int.~J.~Theor.~Phys.}{30}{923-947}{1991}

\ref 
[8] 
C.J.~Isham, ``An introduction to general topology and quantum topology'', 
   in 
    {\it Physics, Geometry and Topology}, 
     proceedings of the 
      Advanced Summer Institute on Physics, Geometry and Topology,  
       held  Banff August 14-26, 1989,
        ed H.C.Lee, pp.129--190 
        (Plenum Press, New York 1990)

\ref 
[9]
R.~Breslav, G.N.~Parfionov and R.R.~Zapatrin,
``Topology measurement within the histories approach'',
\journaldata{Hadronic Journal}{22}{225-239}{1999}
\eprint{quant-ph/9903011}

\ref	
[10]
Ioannis Raptis, ``Algebraic Quantization of Causal Sets'' 
\journaldata{Int. J. Theor. Phys.}{39}{1233-1240}{2000} 
\eprint{gr-qc/9906103}

\ref 
[11] R.D. Sorkin and K.C. Wali (unpublished)

\ref	
[12] 
{\it Encyclopedic Dictionary of Mathematics, Second Edition}
(MIT press, 1993), article 36.D

\ref 
[13]
R.~Geroch, E.H.~Kronheimer and R.~Penrose, ``Ideal Points in Spacetime'',
 \journaldata{Proc. Roy. Soc. Lond.}{A327}{545-567}{1972}

\ref
[14]
Ioannis Raptis and Roman R.~Zapatrin, 
``Quantization of discrete spacetimes and the correspondence principle'',
 \journaldata{Int. J. Theor. Physics}{39}{1}{2000} 
 \eprint{gr-qc/9904079}

\end


(prog1    'now-outlining
  (Outline 
      "
     "......"
      "
   "\\\\message" 
   "\\\\section"
   "\\\\appendi"
   "\\\\Referen"	
   "\\\\Abstrac" 	
      "
   "\\\\subsectio"
   ))